\documentclass{aip-cp}

\usepackage[numbers]{natbib}
\usepackage{rotating}
\usepackage{graphicx}

\usepackage{savesym}
\savesymbol{iint}
\savesymbol{iiint}
\savesymbol{iiiint}
\savesymbol{idotsint}
\usepackage{amsmath}
\usepackage{amssymb}
\restoresymbol{AMS}{iint}
\restoresymbol{AMS}{iiint}
\restoresymbol{AMS}{iiiint}
\restoresymbol{AMS}{idotsint}
\usepackage{multirow}
\usepackage{textcomp}

\begin{document}

\title{A Statistical Method for Corrupt Agents Detection}

\author[aff2]{Yury A. Pichugin}
\eaddress{yury-pichugin@mail.ru}
\author[aff1]{Oleg A. Malafeyev}
\eaddress{malafeyevoa@mail.ru}
\author[aff1]{Denis Rylow\corref{cor1}}

\affil[aff1]{St. Petersburg State University, 7/9 Universitetskaya nab., St. Petersburg, 199034, Russia.}
\affil[aff2]{St. Petersburg State University of Aerospace Instrumentation, 67, Bolshaya Morskaia str, St. Petersburg, Russian Federation}
\corresp[cor1]{Corresponding author: denisrylow@gmail.com}

\maketitle

\begin{abstract} 
The statistical method is used to identify the hidden leaders of the corruption structure. The method is based on principal component analysis (PCA), linear regression, and Shannon information. It is applied to study the time series data of corrupt financial activity. 
 Shannon's quantity of information is specified as a function of two arguments: a vector of hidden corruption factors and a subset of corrupt agents. Several optimization problems are solved to determine the contribution of corresponding corrupt agents to the total illegal behavior. 
An illustrative example is given. A convenient algorithm for computing the covariance matrix with missing data is proposed.

\end{abstract}

\section{INTRODUCTION}

This paper demonstrates new application of principal component analysis (PCA) which was originally proposed by K. Pearson \cite{pearson}. This method is also known as Karhunen-Loeve Transform (see for example \cite{stat_book1,stat_book2}). 
Transition to the independent variables (principal components) allows to reduce the dimension of matrices associated with the problem and to simplify computations when maintaining sufficient precision and relevancy of the results.
This property has made the method particularly popular in time series analysis, where it is known as singular spectrum analysis (ASS). This method has been extensively used in visualization and graphical representation of multidimensional data \cite{Broomhead1,Broomhead2,Ghil1}. In this paper the method is used to formalize a regression model and select relevant data \cite{Muresan,Rao1}. 

The problem of relevant data selection is formalized and solved as an optimization problem of regression analysis. The most well-known criteria of optimality are D-criterion \cite{D_optimal1,D_optimal2}, A-criterion, and G-criterion. A-criterion is iterative one and so quite appealing to the researcher as it allows to easily solve the optimization problem although G.~Seber has criticized its usage \cite{Regression}. In this paper the estimated parameters of regression are principal components that are stochastic in nature. So Shannon's quantity of information in vector form is used as an optimality criterion. 
This approach is similar to D-criterion but, unlike it, is iterative, so the computations can be performed as easily as in the case of A-criterion. Other optimization and relevant numerical techniques could also be potentially used \cite{Tsitouras_optimization_swarm,icnaam_optimiz,kvitko1,kvitko2,simos_icnaam,simos_icnaam2,Kabrits1,zubov1}.

The problem of corruption has been extensively studied mainly in game-theoretic context \cite{Gao2,icnaam_game,Malafeyev_Kolokoltsov_scopus_2010_understanding_game_theory_book,Malafeyev_scopus_2016_parameter_mechanism_design,Malafeyev_scopus_2015_corrupt_inspection,Malafeyev_scopus_2015_interaction_between_anticorruption_authority,Malafeyev_scopus_2017_mean_field_game,Malafeyev_scopus_2016_stochastic_analysis_corrupt_dynamics,Malafeyev_scopus_2014_electric_circutis_analogies_in_economy,Malafeyev_scopus_2015_corrupt_model,Malafeyev_scopus_2016_estimation_of_corruption_indicators}. In this paper the time series of corrupt financial activity data is studied with principal component analysis. The covariant matrix of biased estimators $\bf C$ is calculated with missing data. The matrix $\bf C$ is then used to find the matrix $\bf P$ which is utilized in a linear regression model. 
 Shannon's quantity of information is written as a function of two arguments: a vector of hidden corruption factors and a subset of corrupt agents. Shannon's quantity of information is calculated with biased estimators (matrix $\bf P$) \cite{Yaglom}. Quantity of information maximization problems are formalized for every successive natural number $1,2,3, \dots, m-1,m$. Solutions to these optimization problems allow us to estimate contributions of corrupt agents to the total illegal behavior.  
 A numerical illustrative example is given.

\begin{figure}
\includegraphics{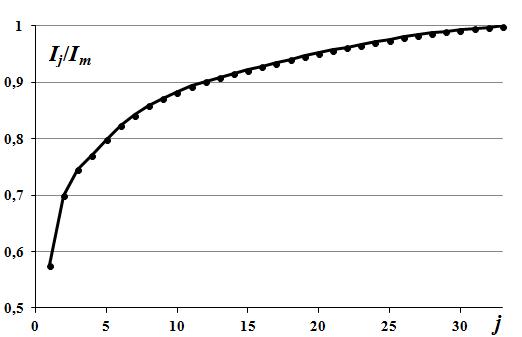}
\caption{Increments in Shannon's quantities of information reflecting contribution of each additional agent to the illegal behavior (hidden corruption factors).}
\label{fig1}
\end{figure} 

\begin{table}
\caption{Monetary values of corrupt financial activity.}
\label{tab1}
\centering
\resizebox{0.95\textwidth}{!}{
\begin{tabular}{|c|c|c|c|c|c|c|c|c|c|c|c|c|c|c|c|c|}
 \hline
\multirow{2}{*}&\multicolumn{7}{|c|}{Periods of observation} & \multirow{2}{*}&\multicolumn{7}{|c|}{Periods of observation}\\ \hline
Agent &1	&2	&3	&4	&5	&6	&7  & Agent &1	&2	&3	&4	&5	&6	&7  \\ \hline
1	&113.1	&50.7	&100.2	&130.5	&129.6	&134.1	&76.2 &
2	&143.4	&96.3	&136.2	&0		&146.4	&163.2	&0\\ 
3	&147	&89.4	&137.4	&158.4	&160.2	&159	&120.9&
4	&110.7	&89.4	&128.4	&48.9	&150.3	&163.2	&112.8\\ 
5	&130.5	&70.2	&121.2	&133.5	&138	&142.5	&75.9&
6	&117.6	&221.1	&121.2	&121.2	&152.1	&0		&100.5\\ 
7	&133.2	&92.4	&126.6	&0		&154.8	&0		&0&
8	&108.9	&0		&121.2	&146.1	&138.3	&150.6	&103.2\\ 
9	&104.7	&96.9	&115.8	&108.6	&146.4	&146.4	&100.5&
10	&126.9	&85.8	&118.2	&143.4	&147.9	&147.9	&92.7\\ 
11	&102.3	&39	9	&3.3		&114.9	&105	&115.2	&64.8&
12	&131.1	&96.6	&125.1	&151.5	&144	&155.1	&109.2\\ 
13	&126.9	&0		&118.5	&0		&138	&155.1	&102.6&
14	&105.6	&76.8	&108.6	&104.7	&130.5	&139.5	&81.6\\ 
15	&0		&93.6	&100.5	&137.7	&141.6	&147	&105.3&
16	&113.4	&88.2	&110.1	&138.9	&0		&147	&102.3\\ 
17	&117.6	&53.1	&104.7	&126.6	&0		&137.4	&72&
18	&0		&0		&117.6	&154.5	&150	&159.6	&102.3\\ 
19	&0		&149.4	&114.9	&158.4	&120.9	&146.4	&105.3&
20	&129.6	&91.8	&110.4	&162.6	&150	&150.6	&92.7\\ 
21	&96.6	&0		&110.4	&139.2	&138.3	&147.9	&92.7&
22	&105	&0		&109.8	&0		&150	&147.9	&96.9\\ 
23	&105.9	&0		&102	&0		&125.1	&135.3	&74.1&
24	&129.6	&91.8	&125.1	&126.9	&149.4	&153.3	&106.5\\ 
25	&118.8	&0		&138.9	&150.6	&0		&159	&110.4&
26	&138.3	&101.1	&137.7	&168.3	&124.8	&163.2	&106.5\\ 
27	&105	&97.8	&0		&150.6	&153.9	&167.1	&109.5&
28	&129.9	&0		&105.6	&146.4	&149.7	&154.8	&96.9\\ 
29	&81.3	&5.1		&71.1	&92.7	&76.5	&92.7	&14.1&
30	&113.4	&0		&1.8		&121.5	&0		&142.5	&80.4\\ 
31	&136.2	&96.3	&25.5	&154.8	&158.1	&159	&113.7&
32	&121.2	&82.8	&0		&151.8	&153.9	&0		&0\\ 
33	&93.3	&42.3	&17.1	&113.1	&87.6	&109.2	&63.9 & 
~	&~	&	&	&	&	&	&\\ \hline
\end{tabular}
}
\end{table}

\section{FORMAL MODEL}

Let $M=\{ 1,2, \dots, m\}$ be an ordered set of $m$ corrupt agents. Divide the observed time interval into equal periods of time of the length $\tau$. Denote the number of these intervals by $J$, so that the intervals are numbered accordingly by $1,2, \dots, J$. 
Each corrupt agent $i$, $i = 1,2,\dots,m$, is engaged in some financial transactions at each time interval $j$, $j=1,2,\dots,J$. Let $x_{ij}$ be a total net value of all financial transaction of the agent $i$ during the time period $j$. As agents' illegal and corrupt behavior is studied, assume that $x_{ij}=0$ if no transactions have signs of being illegal. Let $\bf N$ be an integer matrix such that $n_{ij}=1$ if $x_{ij}\not=0$ and $n_{ij}=0$ otherwise. 

Denote by ${\bf Y } = \{ y_{i,j} \}$ a matrix of deviations from mean values. Let $n_{i,*} = \sum_{j = 1}^{J} n_{i,j}$. If $n_{i,*} \not= 0$, then 
\begin{equation}
\label{one}
y_{i,j} = x_{i,j} - \frac{\sum_{j = 1}^{J} x_{i,j}}{n_{i,*}}.
\end{equation}
If $n_{i,*} =0$ for some $i$, then $y_{i,j} = 0$ for all possible $j$ and this index $i$. Both matrices $\bf Y$ and $\bf N$ are $m \times J$ dimensional. Assume that matrices $\bf N$ and $\bf Y$ do not contain any zero rows. Then it is possible to compute a $m \times m $ matrix ${\bf C} = c_{ij}$ as follows: if $ [ {\bf NN}^T ]_{i,j} \not= 0$, then
\begin{equation}
\label{two}
c_{i,j} =  \frac{[ {\bf YY}^T ]_{i,j} }{ [ {\bf NN}^T ]_{i,j}},
\end{equation}
otherwise $c_{i,j} = 0$. Here $T$ is the transpose operator and $[ \cdot ]_{i,j}$ denotes the element of the corresponding matrix with indices $i,j$. The matrix $\bf C$ obtained in this manner is a matrix of biased covariance estimates of vectors $x_j = (x_{1,j}, x_{2,j}, \dots, x_{m,j})$, $j = 1,2, \dots, J$.  
Let $\bf R $ be an orthogonal $m \times m$ matrix such that $\textbf{R}^T \textbf{C R} = diag(\lambda_1, \lambda_2, \dots, \lambda_m ) = \Lambda $ and $\lambda_1 \geq \lambda_1 \geq \dots \geq \lambda_m$. 

Choose $k \leq m$ values $\lambda_i$. Assume that the relative dispersion $\frac{\sum_{i = 1}^{k} \lambda_i}{\sum_{i = 1}^{m} \lambda_i}$ is large.
Then matrix $\textbf{P}$ is $m \times k$ matrix such that $\textbf{P}^T \textbf{C P} = diag(\lambda_1, \lambda_2, \dots, \lambda_k ) = \Lambda_* $. The matrix $\textbf{P}$ contains the first $k$ columns of $\bf R$.
Let $\bf y$ be some column of the matrix $\bf Y$. Consider a simple regression model 
\begin{equation}
\label{three}
\bf y = P \cdot z + \varepsilon,
\end{equation}
where $\bf z$ is a $k$-dimensional vector and ${\varepsilon}$ is a $m$-dimensional vector of regression errors. The vector $\bf z$ is a well-known vector of principal components. In our case it can be interpreted as a vector of hidden corruption factors.  Assume that components of $ \varepsilon$ do not correlate and the covariance matrix of the vector $ \varepsilon$ is $\sigma^2 I$, where $I$ is the identity matrix.

Note however that biased estimators should be used to compute Shannon's quantity of information. So deviation is calculated as: 
$
\sigma^2 = \frac{tr \textbf{Y}^T ( \textbf{I} - \textbf{PP}^T)	\textbf{Y}}{mJ}.
$
Let $q \subset \{1,2, \dots, m \}$. Denote by $\textbf{y}_q$ a $|q|$-dimensional vector 
such that 
$\textbf{y}_q = (y_l)_{l \in q}$, where $y_l$ is a component of the vector $\bf y$ with index $l$; the order of components is preserved.
Let $\textbf{P}_q$ be a $|q| \times m$-dimensional matrix that is obtained from the matrix $\bf P$ by deleting all rows whose indices are not in $q$.
The quantity of information the vector $\textbf{y}_q$ contains relative to the vector $\bf z$ can be written as:  
\begin{equation}
\label{four}
I(\textbf{y}_q, \textbf{z}) = 0.5 \cdot log_{\alpha} (det( \textbf{I} + \sigma^{-2} \Lambda_*  \textbf{P}_q \textbf{P}_q^T)).
\end{equation}

Matrix $ \sigma^{-2}  \textbf{P}_q \textbf{P}_q^T$ is a Fisher information matrix. So $I(\textbf{y}_q, \textbf{z}) $ can be considered as a function of only one argument $\bf q$. Select several subsets $q$ of vector $\bf y$ components (elements with numbers $\{ 1,2,3, \dots, m \}$) by incrementally increasing $|q|$ from $0$ to $m$. Maximize $I(\textbf{y}_q, \textbf{z})$ for each such subset. The result is a sequence of values  $I(\textbf{y}_q, \textbf{z})$ , i.~e. $I_1 \leq I_2 \leq \dots \leq I_m$. Values $\Delta I_1 = I_1$,  $\Delta I_2 = I_2 - I_1$, ..., $\Delta I_m = I_m - I_{m-1}$ are increments in Shannon's quantity of information. Note that $\Delta I_1 \geq \Delta I_2 \geq \dots \geq \Delta I_m $. 
The logarithm base $\alpha $ in (\ref{four}) can be chosen arbitrarily if relative values $\frac{I_j}{I_m}$ and $\frac{\Delta I_j}{I_m}$ are considered.
Here $I_m = I(\textbf{y},\textbf{z})$ is the full quantity of information in the model (\ref{three}).

If $k<<m$ and $j = j_*$ is relatively small, then the fraction $\frac{I_j}{I_m}$ can become quite large (for example $0.8$--$0.9$) or the value $\frac{\Delta I_j}{I_m}$ can become quite small (for example $0.1$). The subset $\bf q$ with the corresponding $j_* = |\textbf{q}|$ is the subset of corrupt agents that constitute the core of the corruption structure. The choice of the subset $q$ is equivalent to the choice of agents that should be 
closely monitored for possible corrupt behavior.

{\bf Remark.} The requirement of absence of zero rows in matrices $\bf Y$ and $\bf N$ can be relaxed. Presence of zero rows in these matrices leads to the presence of zero rows in $\bf P$, but the indices of these rows and the corresponding agents are the last in the successive data selection process.

\section{NUMERICAL EXAMPLE}

Consider an illustrative example. Table \ref{tab1} contains observed values ($x_{ij}$). The number of agents $m=33$, the number of time periods $J=7$. By replacing all non-zero values with $1$ we get the matrix $\bf N$. Then matrices $\bf Y$ and $\bf C$ can be found using (\ref{one}) and (\ref{two}). Note that $\sum_{i=1}^{m} \lambda_i = tr \textbf{C}$. So the relative dispersion of the $k$ components is $\frac{\sum_{i=1}^{k} \lambda_i}{tr \textbf{C}}$. 
We use von Mises iteration to calculate the values $\lambda_i$ and columns of the matrix $\bf P$.
For a example, if $k = 3$, then  $\sum_{i=1}^{k} \lambda_i / \sum_{i=1}^{m} \lambda_i = 0.98$ and $tr( \textbf{C} )^{-1} \Lambda_* = diag(0.65, 0.21, 0.12)$. So the matrix $\bf P$ is of size  $33 \times 3$.

The figure \ref{fig1} contains the sequence of values $\frac{I_j}{I_m}$. The first third of all agents gives almost $90\%$ of all information about hidden corruption factors. Indices of those agents are $30, 31, 33, 22, 6, 23, 18, 4, 28, 21, 27$. Unexpectedly, agent $30$ contributes the most to the corrupt behavior. This is completely counter intuitive if only raw data in the table \ref{tab1} is considered.

\section{CONCLUSION}
Method of statistical analysis is applied to the problem of corrupt behavior detection. Principal component analysis, linear regression, Shannon's information are used to identify the hidden leaders of the corruption structure. Shannon's quantities of information are maximized to determine the contribution of corrupt agents to the total illegal behavior. A numerical example is provided. A convenient algorithm for computation of covariant matrix is proposed.


\nocite{*}

\hyphenpenalty=10000
\bibliography{malafeyev_scopus_web_of_science(Science_Index)}%

\begin{thebibliography}{100}
\def\selectlanguageifdefined#1{
\expandafter\ifx\csname date#1\endcsname\relax
\else\selectlanguage{#1}\fi}
\providecommand*{\href}[2]{{\small #2}}
\providecommand*{\url}[1]{{\small #1}}
\providecommand*{\BibUrl}[1]{\url{#1}}
\providecommand{\BibAnnote}[1]{}
\providecommand*{\BibEmph}[1]{#1}
\ProvideTextCommandDefault{\cyrdash}{\iflanguage{russian}{\hbox
  to.8em{--\hss--}}{\textemdash}}
\providecommand*{\BibDash}{\ifdim\lastskip>0pt\unskip\nobreak\hskip.2em plus
  0.1em\fi
\cyrdash\hskip.2em plus 0.1em\ignorespaces}
\renewcommand{\newblock}{\ignorespaces}

\bibitem{pearson}
\selectlanguageifdefined{english}
\BibEmph{Pearson~K.} On lines and planes of closest fit to systems of points in
  space~// \BibEmph{Philosophical Magazine}. \BibDash
\newblock 1901. \BibDash
\newblock Vol.~2. \BibDash
\newblock P.~559--572.

\bibitem{stat_book1}
\selectlanguageifdefined{english}
\BibEmph{Karhunen~Kari}. Uber lineare Methoden in der
  Wahrscheinlichkeitsrechnung~// \BibEmph{Annales Academiae scientiarum
  Fennicae. Series A. 1, Mathematica-physica}. \BibDash
\newblock 1947. \BibDash
\newblock Vol.~37. \BibDash
\newblock P.~1--79.

\bibitem{stat_book2}
\selectlanguageifdefined{english}
\BibEmph{Loeve~M.} Probability Theory. \BibDash
\newblock 4 edition. \BibDash
\newblock Springer-Verlag New York, 1978. \BibDash
\newblock Vol.~2 of \BibEmph{Graduate Texts in Mathematics}.

\bibitem{Broomhead1}
\selectlanguageifdefined{english}
\BibEmph{Broomhead~D.S., King~G.P.} Extracting qualitative dynamics from
  experimental data~// \BibEmph{Physica D}. \BibDash
\newblock 1986. \BibDash
\newblock Vol.~20. \BibDash
\newblock P.~21--236.

\bibitem{Broomhead2}
\selectlanguageifdefined{english}
\BibEmph{Broomhead~D.S., King~G.P.} On the qualitative analysis of experimental
  dynamical systems~// Nonlinear Phenomena and Chaos~/ Ed.\ by\ S.~Sarkar.
  \BibDash
\newblock Adam Hilger, Bristol, 1986. \BibDash
\newblock Vol.~20. \BibDash
\newblock P.~113--144.

\bibitem{Ghil1}
\selectlanguageifdefined{english}
\BibEmph{Ghil~M., Vautard~R.} Interdecadal oscillations and the warming trend
  in global temperature time series~// \BibEmph{Nature}. \BibDash
\newblock 1991. \BibDash
\newblock Vol. 350. \BibDash
\newblock P.~324--327.

\bibitem{Muresan}
\selectlanguageifdefined{english}
\BibEmph{Muresan~D.~D., Parks~T.~W.} Adaptive principal components and image
  denoising~// Proceedings 2003 International Conference on Image Processing
  (Cat. No.03CH37429). \BibDash
\newblock Vol.~1. \BibDash
\newblock 2003. \BibDash
\newblock P.~I--101--4 vol.1.

\bibitem{Rao1}
\selectlanguageifdefined{english}
\BibEmph{Rao~Kamisetty~Ramam, Yip~Patrick~C.} The Transform and Data
  Compression Handbook. \BibDash
\newblock CRC Press, 2000.

\bibitem{D_optimal1}
\selectlanguageifdefined{english}
\BibEmph{Mitchell~Toby~J.} An Algorithm for the Construction of "D-Optimal"
  Experimental Designs~// \BibEmph{Technometrics}. \BibDash
\newblock 1974. \BibDash
\newblock Vol.~16, no.~2. \BibDash
\newblock P.~203--210.

\bibitem{D_optimal2}
\selectlanguageifdefined{english}
\BibEmph{Mitchell~Toby~J.} Computer Construction of "D-Optimal" First-Order
  Designs~// \BibEmph{Technometrics}. \BibDash
\newblock 1974. \BibDash
\newblock Vol.~16, no.~2. \BibDash
\newblock P.~211--220.

\bibitem{Regression}
\selectlanguageifdefined{english}
\BibEmph{Seber~George A.~F., Lee~Alan~J.} Linear Regression Analysis. \BibDash
\newblock 2 edition. \BibDash
\newblock John Wiley \& Sons, 2003.

\bibitem{Tsitouras_optimization_swarm}
\selectlanguageifdefined{english}
\BibEmph{Alexandridis~Alex, Famelis~Ioannis~Th., Tsitouras~Charalambos}.
  Particle swarm optimization for complex nonlinear optimization problems~//
  \BibEmph{AIP Conference Proceedings}. \BibDash
\newblock 2016. \BibDash
\newblock Vol. 1738, no. 480120.

\bibitem{icnaam_optimiz}
\selectlanguageifdefined{english}
\BibEmph{e~Silva~Eliana~Costa, Correia~Aldina, Lopes~Isabel~Cristina}.
  Optimization in generalized linear models: A case study~// \BibEmph{AIP
  Conference Proceedings}. \BibDash
\newblock 2016. \BibDash
\newblock Vol. 1738, no. 300002.

\bibitem{kvitko1}
\selectlanguageifdefined{english}
\BibEmph{Kvitko~Alexander}. A method for solving a local boundary problem for
  nonlinear controlled system~// \BibEmph{AIP Conference Proceedings}. \BibDash
\newblock 2015. \BibDash
\newblock Vol. 1648, no. 450002.

\bibitem{kvitko2}
\selectlanguageifdefined{english}
\BibEmph{Kvitko~Alexander}. Syntheses of terminal control for nonlinear
  stationary controlled system under incomplete information~// \BibEmph{AIP
  Conference Proceedings}. \BibDash
\newblock 2016. \BibDash
\newblock Vol. 1738, no. 160002.

\bibitem{simos_icnaam}
\selectlanguageifdefined{english}
An optimized two-step hybrid block method for solving general second order
  initial-value problems of the form y″ = f (x, y, y′)~/ Higinio~Ramos,
  Z.~Kalogiratou, Th.~Monovasilis, T.E.~Simos~// \BibEmph{AIP Conference
  Proceedings}. \BibDash
\newblock 2015. \BibDash
\newblock Vol. 1648, no. 810006.

\bibitem{simos_icnaam2}
\selectlanguageifdefined{english}
A trigonometrically fitted optimized two-step hybrid block method for solving
  initial-value problems of the form y″ = f (x, y, y′) with oscillatory
  solutions~/ Higinio~Ramos, Z.~Kalogiratou, Th.~Monovasilis, T.E.~Simos~//
  \BibEmph{AIP Conference Proceedings}. \BibDash
\newblock 2015. \BibDash
\newblock Vol. 1648, no. 810007.

\bibitem{Kabrits1}
\selectlanguageifdefined{english}
\BibEmph{Kabrits~S.A., Kolpak~E.P.} Finding bifurcation branches in nonlinear
  problems of statics of shells numerically~// 2015 International Conference on
  "Stability and Control Processes" in Memory of V.I. Zubov, SCP 2015 -
  Proceedings. \BibDash
\newblock 2015. \BibDash
\newblock P.~389--391.

\bibitem{zubov1}
\selectlanguageifdefined{english}
Application in practice and optimization of industrial information systems~/
  L.A.~Bondarenko, A.V.~Zubov, V.B.~Orlov et~al.~// \BibEmph{Journal of
  Theoretical and Applied Information Technology}. \BibDash
\newblock 2016. \BibDash
\newblock Vol.~85, no.~3. \BibDash
\newblock P.~305--308.

\bibitem{Gao2}
\selectlanguageifdefined{english}
Cooperation in two-stage games on undirected networks~/ H.~Gao, L.~Petrosyan,
  H.~Qiao, A.~Sedakov~// \BibEmph{Journal of Systems Science and Complexity}.
  \BibDash
\newblock 2017. \BibDash
\newblock Vol.~30, no.~3. \BibDash
\newblock P.~680--693.

\bibitem{icnaam_game}
\selectlanguageifdefined{english}
\BibEmph{Meirong~Wu, Shaochen~Cao, Huazhen~Zhu}. On axiomatizations of the
  Shapley value for bi-cooperative games~// \BibEmph{AIP Conference
  Proceedings}. \BibDash
\newblock 2016. \BibDash
\newblock Vol. 1738, no. 080002.

\bibitem{Malafeyev_Kolokoltsov_scopus_2010_understanding_game_theory_book}
\selectlanguageifdefined{english}
\BibEmph{Malafeyev~O.A., Kolokoltsov~V.N.} Understanding game theory:
  Introduction to the analysis of many agent systems with competition and
  cooperation. \BibDash
\newblock New Jersey~: World Scientific Publishing Co., 2010.

\bibitem{Malafeyev_scopus_2016_parameter_mechanism_design}
\selectlanguageifdefined{english}
\BibEmph{Malafeyev~O., Pichugin~Y., Alferov~G.} Parameters estimation in
  mechanism design~// \BibEmph{Contemporary Engineering Sciences}. \BibDash
\newblock 2016. \BibDash
\newblock Vol.~9, no. 1-4. \BibDash
\newblock P.~175--185.

\bibitem{Malafeyev_scopus_2015_corrupt_inspection}
\selectlanguageifdefined{english}
\BibEmph{Malafeyev~O.A., Alferov~G.V., Maltseva~A.S.} Game-theoretic model of
  inspection by anti-corruption group~// \BibEmph{AIP Conference Proceedings}.
  \BibDash
\newblock 2015. \BibDash
\newblock Vol. 1648, no. 450009.

\bibitem{Malafeyev_scopus_2015_interaction_between_anticorruption_authority}
\selectlanguageifdefined{english}
\BibEmph{Neverova~E.G., Malafeyef~O.A.} A model of interaction between
  anticorruption authority and corruption groups~// \BibEmph{AIP Conference
  Proceedings}. \BibDash
\newblock 2015. \BibDash
\newblock Vol. 1648, no. 450012.

\bibitem{Malafeyev_scopus_2017_mean_field_game}
\selectlanguageifdefined{english}
\BibEmph{Kolokoltsov~V.N., Malafeyev~O.A.}~// \BibEmph{Dynamic Games and
  Applications}. \BibDash
\newblock 2017. \BibDash
\newblock Vol.~7, no.~1. \BibDash
\newblock P.~34--47.

\bibitem{Malafeyev_scopus_2016_stochastic_analysis_corrupt_dynamics}
\selectlanguageifdefined{english}
\BibEmph{Malafeyev~O.A., Redinskikh~N.D.} Stochastic analysis of the dynamics
  of corrupt hybrid networks~// Proceedings of 2016 International Conference
  "Stability and Oscillations of Nonlinear Control Systems" (Pyatnitskiy's
  Conference), STAB 2016. \BibDash
\newblock No. 123354. \BibDash
\newblock New Jersey~: Institute of Electrical and Electronics Engineers Inc.,
  2016. \BibDash
\newblock P.~1--4.

\bibitem{Malafeyev_scopus_2014_electric_circutis_analogies_in_economy}
\selectlanguageifdefined{english}
\BibEmph{Malafeyev~O.A., Redinskikh~N.D., Alferov~G.V.} Electric circuits
  analogies in economics modeling: Corruption networks~// 2014 2nd
  International Conference on Emission Electronics, ICEE 2014. \BibDash
\newblock New Jersey~: Institute of Electrical and Electronics Engineers Inc.,
  2014. \BibDash
\newblock P.~28--32.

\bibitem{Malafeyev_scopus_2015_corrupt_model}
\selectlanguageifdefined{english}
Model of interaction between anticorruption authorities and corruption groups~/
  O.A.~Malafeyev, E.G.~Neverova, G.V.~Alferov, T.E.~Smirnova~// 2015
  International Conference on ``Stability and Control Processes'' in Memory of
  V.I. Zubov, SCP 2015 - Proceedings. \BibDash
\newblock New Jersey~: Institute of Electrical and Electronics Engineers Inc.,
  2015. \BibDash
\newblock P.~488--490.

\bibitem{Malafeyev_scopus_2016_estimation_of_corruption_indicators}
\selectlanguageifdefined{english}
\BibEmph{Pichugin~Y.A., Malafeyev~O.A.} Statistical estimation of corruption
  indicators in the firm~// \BibEmph{Applied Mathematical Sciences}. \BibDash
\newblock 2016. \BibDash
\newblock Vol.~10, no. 41-44. \BibDash
\newblock P.~2065--2073.

\bibitem{Yaglom}
\selectlanguageifdefined{english}
\BibEmph{Gelyofand~I.~M., Yaglom~A.~M.} Computation of the amount of
  information about a stochastic function contained in another such function~//
  \BibEmph{Uspehi Mat. Nauk (N.S.)}. \BibDash
\newblock 1957. \BibDash
\newblock Vol.~12, no. 1(73). \BibDash
\newblock P.~3–52.

\bibitem{Malafeyev_rinz_Dynamic_competitive_systems_1}
\selectlanguageifdefined{english}
\BibEmph{Malafeev~O.~A., Kolokol'czov~V.~N.} Dinamicheskie konkurentny'e
  sistemy' mnogoagentnogo vzaimodejstviya i ix asimptoticheskoe povedenie
  (chast' I)~// \BibEmph{Vestnik grazhdanskix inzhenerov}. \BibDash
\newblock 2010. \BibDash
\newblock no.~4. \BibDash
\newblock P.~144--153.

\bibitem{Malafeyev_rinz_Dynamic_competitive_systems_2}
\selectlanguageifdefined{english}
\BibEmph{Malafeev~O.~A., Kolokol'czov~V.~N.} Dinamicheskie konkurentny'e
  sistemy' mnogoagentnogo vzaimodejstviya i ix asimptoticheskoe povedenie
  (chast' II)~// \BibEmph{Vestnik grazhdanskix inzhenerov}. \BibDash
\newblock 2011. \BibDash
\newblock no.~1. \BibDash
\newblock P.~134--145.

\bibitem{Malafeyev_rinz_control_process_three_agents}
\selectlanguageifdefined{english}
\BibEmph{Malafeev~O.~A., Sosnina~V.~V.} Model' upravleniya proczessom
  kooperativnogo trexagentnogo vzaimodejstviya~// \BibEmph{Problemy' mexaniki i
  upravleniya: Nelinejny'e dinamicheskie sistemy'}. \BibDash
\newblock 2007. \BibDash
\newblock no.~39. \BibDash
\newblock P.~131--144.

\bibitem{Malafeyev_rinz_math_computer_modelling_many_agents}
\selectlanguageifdefined{english}
\BibEmph{Malafeev~O.~A., Zubova~A.~F.} Matematicheskoe i komp'yuternoe
  modelirovanie soczial'no-e'konomicheskix sistem na urovne mnogoagentnogo
  vzaimodejstviya (vvedenie v problemy' ravnovesiya, ustojchivosti,
  nadezhnosti). \BibDash
\newblock Sankt-Peterburg~: Mobil'nost'-plyus, 2006.

\bibitem{Malafeyev_rinz_postman_problem}
\selectlanguageifdefined{english}
\BibEmph{Malafeev~O.~A., Grigor'eva~K.~V.} Dinamicheskij proczess
  kooperativnogo vzaimodejstviya v mnogokriterial'noj (mnogoagentnoj) zadache
  pochtal'ona~// \BibEmph{Vestnik grazhdanskix inzhenerov}. \BibDash
\newblock 2011. \BibDash
\newblock no.~1. \BibDash
\newblock P.~150--156.

\bibitem{Malafeyev_rinz_conflict_systems}
\selectlanguageifdefined{english}
\BibEmph{Malafeev~O.~A.} Upravlyaemy'e konfliktny'e sistemy'. \BibDash
\newblock Sankt-Peterburg~: SPbGU, 2000.

\bibitem{Malafeyev_rinz_venture_building}
\selectlanguageifdefined{english}
\BibEmph{Malafeev~O.~A., Zenovich~O.~S., Sevek~V.~K.} Mnogoagentnoe
  vzaimodejstvie v dinamicheskoj zadache upravleniya venchurny'mi
  stroitel'ny'mi proektami~// \BibEmph{E'konomicheskoe vozrozhdenie Rossii}.
  \BibDash
\newblock 2012. \BibDash
\newblock no.~1. \BibDash
\newblock P.~124--131.

\bibitem{Malafeyev_rinz_efficiency_city_building}
\selectlanguageifdefined{english}
\BibEmph{Malafeev~O.~A., Drozdova~I.~V., Parshina~L.~G.} E'ffektivnost'
  variantov rekonstrukczii gorodskoj zhiloj zastrojki~//
  \BibEmph{E'konomicheskoe vozrozhdenie Rossii}. \BibDash
\newblock 2008. \BibDash
\newblock no.~3. \BibDash
\newblock P.~63--67.

\bibitem{Malafeyev_rinz_investing_competition}
\selectlanguageifdefined{english}
\BibEmph{Malafeev~O.~A., Paxar~O.~V.} Dinamicheskaya nestaczionarnaya zadacha
  investirovaniya proektov v usloviyax konkurenczii~// \BibEmph{Problemy'
  mexaniki i upravleniya: Nelinejny'e dinamicheskie sistemy'}. \BibDash
\newblock 2009. \BibDash
\newblock no.~41. \BibDash
\newblock P.~103--108.

\bibitem{Malafeyev_rinz_investing_influence_factors_innovate}
\selectlanguageifdefined{english}
\BibEmph{Malafeev~O.~A., Gordeev~D.~A., Titova~N.~D.} Probabilistic and
  deterministic model of the influence factors on the activities of the
  organization to innovate~// \BibEmph{E'konomicheskoe vozrozhdenie Rossii}.
  \BibDash
\newblock 2011. \BibDash
\newblock no.~1. \BibDash
\newblock P.~73--82.

\bibitem{Malafeyev_rinz_coalition_investment_project}
\selectlanguageifdefined{english}
\BibEmph{Malafeev~O.~A., Grigor'eva~K.~V., Ivanov~A.~S.} Statisticheskaya
  koaliczionnaya model' investirovaniya innovaczionny'x proektov~//
  \BibEmph{E'konomicheskoe vozrozhdenie Rossii}. \BibDash
\newblock 2011. \BibDash
\newblock no.~4. \BibDash
\newblock P.~90--98.

\bibitem{Malafeyev_rinz_company_development}
\selectlanguageifdefined{english}
\BibEmph{Malafeev~O.~A., CHerny'x~K.~S.} Matematicheskoe modelirovanie
  razvitiya kompanii~// \BibEmph{E'konomicheskoe vozrozhdenie Rossii}. \BibDash
\newblock 2004. \BibDash
\newblock no.~1. \BibDash
\newblock P.~60--66.

\bibitem{Malafeyev_rinz_stochastic_innovative_product}
\selectlanguageifdefined{english}
\BibEmph{Malafeev~O.~A., Gordeev~D.~A., Titova~N.~D.} Stoxasticheskaya model'
  prinyatiya resheniya o vy'vode na ry'nok innovaczionnogo produkta~//
  \BibEmph{Vestnik grazhdanskix inzhenerov}. \BibDash
\newblock 2011. \BibDash
\newblock no.~2. \BibDash
\newblock P.~161--166.

\bibitem{Malafeyev_rinz_game_theory_for_everyone}
\selectlanguageifdefined{english}
\BibEmph{Malafeev~O.~A., Kolokol'czov~V.~N.} Matematicheskoe modelirovanie
  mnogoagentny'x sistem konkurenczii i kooperaczii (teoriya igr dlya vsex).
  \BibDash
\newblock Sankt-Peterburg~: Lan', 2012.

\bibitem{Malafeyev_rinz_competition_many_agents_auction}
\selectlanguageifdefined{english}
\BibEmph{Malafeev~O.~A., Griczaj~K.~N.} Zadacha konkurentnogo upravleniya v
  modeli mnogoagentnogo vzaimodejstviya aukczionnogo tipa~// \BibEmph{Problemy'
  mexaniki i upravleniya: Nelinejny'e dinamicheskie sistemy'}. \BibDash
\newblock 2007. \BibDash
\newblock no.~39. \BibDash
\newblock P.~36--45.

\bibitem{Malafeyev_rinz_problems_reconstruction_city}
\selectlanguageifdefined{english}
\BibEmph{Malafeev~O.~A., Akulenkova~I.~V., Drozdov~G.~D.} Problemy'
  rekonstrukczii zhilishhno-kommunal'nogo xozyajstva megapolisa. \BibDash
\newblock Sankt-Peterburg~: Sankt-Peterburgskij gosudarstvenny'j universitet
  servisa i e'konomiki, 2007.

\bibitem{Malafeyev_rinz_equilibrium_control_networks}
\selectlanguageifdefined{english}
\BibEmph{Malafeev~O.~A., Parfenov~A.~P.} Ravnovesie i kompromissnoe upravlenie
  v setevy'x modelyax mnogoagentnogo vzaimodejstviya~// \BibEmph{Problemy'
  mexaniki i upravleniya: Nelinejny'e dinamicheskie sistemy'}. \BibDash
\newblock 2007. \BibDash
\newblock no.~39. \BibDash
\newblock P.~154--167.

\bibitem{Malafeyev_rinz_hamilton_jacobi_solution}
\selectlanguageifdefined{english}
\BibEmph{Malafeyev~O.~A., Troeva~M.~S.} A weak solution of Hamilton-Jacobi
  equation for a differential two-person zero-sum game~// Preprints of the
  Eighth International Symposium on Differential Games and Applications.
  \BibDash
\newblock Maastricht~: Rijksuniversiteit te Utrecht, Universiteit Maastricht,
  Rijksuniversiteit te Groningen, 1998. \BibDash
\newblock P.~366--369.

\bibitem{Malafeyev_rinz_reconstruction_city_competition}
\selectlanguageifdefined{english}
\BibEmph{Malafeev~O.~A., Drozdova~I.~V., Drozdov~G.~D.} Modelirovanie
  proczessov rekonstrukczii zhilishhno-kommunal'nogo xozyajstva megapolisa v
  usloviyax konkurentnoj sredy'. \BibDash
\newblock Sankt-Peterburg~: Sankt-Peterburgskij gosudarstvenny'j universitet
  servisa i e'konomiki, 2008.

\bibitem{Malafeyev_rinz_many_agents_network_corruption}
\selectlanguageifdefined{english}
Kompromiss i ravnovesie v modelyax mnogoagentnogo upravleniya v korrupczionnoj
  seti socziuma~/ O.~A.~Malafeev, D.~S.~Bojczov, N.~D.~Redinskix,
  E.~G.~Neverova~// \BibEmph{Molodoj ucheny'j}. \BibDash
\newblock 2014. \BibDash
\newblock no. 10 (69). \BibDash
\newblock P.~14--17.

\bibitem{Malafeyev_rinz_american_english_corruption}
\selectlanguageifdefined{english}
\BibEmph{Malafeev~O.~A., Ry'lov~D.~S.} Formalizacziya grazhdanskogo proczessa
  pri anglijskom i amerikanskom sudebnom pravile v vide bajesovskoj igry'~//
  \BibEmph{Molodoj ucheny'j}. \BibDash
\newblock 2016. \BibDash
\newblock no.~13. \BibDash
\newblock P.~46--50.

\bibitem{Malafeyev_rinz_competition_auction}
\selectlanguageifdefined{english}
\BibEmph{Malafeev~O.~A., Griczaj~K.~N.} Konkurentnoe upravlenie v modelyax
  aukczionov~// \BibEmph{Problemy' mexaniki i upravleniya: Nelinejny'e
  dinamicheskie sistemy'}. \BibDash
\newblock 2004. \BibDash
\newblock no.~36. \BibDash
\newblock P.~74--82.

\bibitem{Malafeyev_rinz_conflict_market_entry}
\selectlanguageifdefined{english}
\BibEmph{Malafeev~O.~A., Ershova~T.~A.} Konfliktny'e upravleniya v modeli
  vxozhdeniya v ry'nok~// \BibEmph{Problemy' mexaniki i upravleniya:
  Nelinejny'e dinamicheskie sistemy'}. \BibDash
\newblock 2004. \BibDash
\newblock no.~36. \BibDash
\newblock P.~19--27.

\bibitem{Malafeyev_rinz_methods_postman}
\selectlanguageifdefined{english}
\BibEmph{Malafeev~O.~A., Grigor'eva~K.~V.} Metody' resheniya dinamicheskoj
  mnogokriterial'noj zadachi pochtal'ona~// \BibEmph{Vestnik grazhdanskix
  inzhenerov}. \BibDash
\newblock 2011. \BibDash
\newblock no.~4. \BibDash
\newblock P.~156--161.

\bibitem{Malafeyev_rinz_equilibrium_methods_conflict_control_system}
\selectlanguageifdefined{english}
\BibEmph{Malafeev~O.~A., Troeva~M.~S.} Ustojchivost' i nekotory'e chislenny'e
  metody' v konfliktno upravlyaemy'x sistemax. \BibDash
\newblock YAkutsk~: YAkutskij gosudarstvenny'j universitet im. M. K. Ammosova,
  1998.

\bibitem{Malafeyev_rinz_agricultural_production}
\selectlanguageifdefined{english}
Matematicheskoe modelirovanie proczessov v agropromy'shlennom proizvodstve~/
  O.~A.~Malafeev, V.~S.~SHkrabak, A.~V.~Skrobach, V.~F.~Skrobach. \BibDash
\newblock Sankt-Peterburg~: Sankt-Peterburgskij gosudarstvenny'j agrarny'j
  universitet, 2000.

\bibitem{Malafeyev_rinz_existence_general_value}
\selectlanguageifdefined{english}
\BibEmph{Malafeev~O.~A.} O sushhestvovanii obobshhennogo znacheniya
  dinamicheskoj igry'~// \BibEmph{Vestnik Sankt-Peterburgskogo universiteta.
  Seriya 1. Matematika. Mexanika. Astronomiya}. \BibDash
\newblock 1972. \BibDash
\newblock no.~4. \BibDash
\newblock P.~41--46.

\bibitem{Malafeyev_rinz_conflict_situation_solution_2}
\selectlanguageifdefined{english}
\BibEmph{Malafeev~O.~A., Murav'ev~A.~I.} Matematicheskie modeli konfliktny'x
  situaczij i ix razreshenie. \BibDash
\newblock Sankt-Peterburg~: Sankt-Peterburgskij gosudarstvenny'j universitet
  e'konomiki i finansov, 2001. \BibDash
\newblock Vol.~2.

\bibitem{Malafeyev_rinz_conflict_situation_solution_1}
\selectlanguageifdefined{english}
\BibEmph{Malafeev~O.~A., Murav'ev~A.~I.} Matematicheskie modeli konfliktny'x
  situaczij i ix razreshenie. \BibDash
\newblock Sankt-Peterburg~: Sankt-Peterburgskij gosudarstvenny'j universitet
  e'konomiki i finansov, 2000. \BibDash
\newblock Vol.~1.

\bibitem{Malafeyev_rinz_control_city_building}
\selectlanguageifdefined{english}
\BibEmph{Malafeev~O.~A., Drozdov~G.~D.} Modelirovanie proczessov v sisteme
  upravleniya gorodskim stroitel'stvom. \BibDash
\newblock Sankt-Peterburg~: Sankt-Peterburgskij gosudarstvenny'j
  arxitekturno-stroitel'ny'j universitet, 2001. \BibDash
\newblock Vol.~1.

\bibitem{Malafeyev_rinz_corruption_contracts}
\selectlanguageifdefined{english}
\BibEmph{Malafeev~O.~A., Koroleva~O.~A.} Model' korrupczii pri zaklyuchenii
  kontraktov~// Proczessy' upravleniya i ustojchivost' Trudy' XXXIX
  mezhdunarodnoj nauchnoj konferenczii aspirantov i studentov pod redakcziej N.
  V. Smirnova, G. SH. Tamasyana. \BibDash
\newblock Sankt-Peterburg~: Sankt-Peterburgskij gosudarstvenny'j universitet,
  2008. \BibDash
\newblock P.~446--449.

\bibitem{Malafeyev_rinz_conflict_situation_socio_economic_system}
\selectlanguageifdefined{english}
\BibEmph{Malafeev~O.~A., Murav'ev~A.~I.} Modelirovanie konfliktny'x situaczij v
  soczial'no-e'konomicheskix sistemax. \BibDash
\newblock Sankt-Peterburg~: Sankt-Peterburgskij gosudarstvenny'j universitet
  e'konomiki i finansov, 1998.

\bibitem{Malafeyev_rinz_many_agents_insurance}
\selectlanguageifdefined{english}
\BibEmph{Malafeev~O.~A., Drozdov~G.~D.} Modelirovanie mnogoagentnogo
  vzaimodejstviya proczessov straxovaniya. \BibDash
\newblock Sankt-Peterburg~: Sankt-Peterburgskij gosudarstvenny'j universitet
  servisa i e'konomiki, 2010.

\bibitem{Malafeyev_rinz_equilibrium_oscillation_economic_model}
\selectlanguageifdefined{english}
\BibEmph{Malafeev~O.~A., Zubova~A.~F.} Ustojchivost' po Lyapunovu i
  kolebatel'nost' v e'konomicheskix modelyax. \BibDash
\newblock Sankt-Peterburg~: Sankt-Peterburgskij gosudarstvenny'j universitet,
  2001.

\bibitem{Malafeyev_rinz_strategy_n_game}
\selectlanguageifdefined{english}
\BibEmph{Malafeev~O.~A., Bure~V.~M.} Soglasovannaya strategiya v
  povtoryayushhixsya konechny'x igrax N licz~// \BibEmph{Vestnik
  Sankt-Peterburgskogo universiteta. Seriya 1. Matematika. Mexanika.
  Astronomiya}. \BibDash
\newblock 1995. \BibDash
\newblock no.~1. \BibDash
\newblock P.~120--122.

\bibitem{Malafeyev_rinz_value_existence_pursuit}
\selectlanguageifdefined{english}
\BibEmph{Malafeev~O.~A.} O sushhestvovanii znacheniya igry' presledovaniya~//
  \BibEmph{Sibirskij zhurnal issledovaniya operaczij}. \BibDash
\newblock 1970. \BibDash
\newblock no.~5. \BibDash
\newblock P.~25--36.

\bibitem{Malafeyev_rinz_phd_thesis}
\selectlanguageifdefined{english}
\BibEmph{Malafeev~O.~A.} Konfliktno upravlyaemy'e proczessy' so mnogimi
  uchastnikami~: Disssertacziya na soiskanie uchenoj stepeni doktora
  fiziko-matematicheskix nauk~/ O.~A.~Malafeev~; Leningradskij gosudarstvenny'j
  universitet. \BibDash
\newblock Leningrad, 1987.

\bibitem{Malafeyev_rinz_optimization_controllable_dynamic_process}
\selectlanguageifdefined{english}
\BibEmph{Malafeev~O.~A.} Ustojchivost' reshenij zadach mnogokriterial'noj
  optimizaczii i konfliktno upravlyaemy'e dinamicheskie proczessy'. \BibDash
\newblock Sankt-Peterburg~: Sankt-Peterburgskij gosudarstvenny'j universitet,
  1990.

\bibitem{Malafeyev_rinz_corruption_auction_first_price}
\selectlanguageifdefined{english}
Korrupcziya v modelyax aukcziona pervoj czeny'~/ O.~A.~Malafeev,
  N.~D.~Redinskix, G.~V.~Alfyorov, T.~E.~Smirnova~// Upravlenie v morskix i
  ae'rokosmicheskix sistemax (UMAS-2014) 7-ya Rossijskaya mul'tikonferencziya
  po problemam upravleniya: materialy' konferenczii. GNCZ RF OAO
  ``CZentral'ny'j nauchno-issledovatel'skij institut ``E'lektropribor''.
  \BibDash
\newblock Sankt-Peterburg~: Konczern ``CZentral'ny'j nauchno-issledovatel'skij
  institut ``E'lektropribor'', 2014. \BibDash
\newblock P.~141--146.

\bibitem{Malafeyev_rinz_network_investment_corruption}
\selectlanguageifdefined{english}
\BibEmph{Malafeev~O.~A., Redinskix~N.~D., Smirnova~T.~E.} Setevaya model'
  investirovaniya proektov s korrupcziej~// Proczessy' upravleniya i
  ustojchivost' Trudy' XLVI mezhdunarodnoj nauchnoj konferenczii aspirantov i
  studentov. \BibDash
\newblock Sankt-Peterburg~: Sankt-Peterburgskij gosudarstvenny'j universitet,
  2015. \BibDash
\newblock P.~659--664.

\bibitem{Malafeyev_rinz_firm_bankruptcy}
\selectlanguageifdefined{english}
\BibEmph{Malafeev~O.~A., Pichugin~YU.~A.} Ob oczenke riska bankrotstva
  firmy'~// Tezisy' dokladov VI Mezhdunar. konf. ``Dinamicheskie sistemy':
  ustojchivost', upravlenie, optimizacziya'' (DSSCO’13). \BibDash
\newblock Minsk~: Belorusskij gosudarstvenny'j universitet, 2013. \BibDash
\newblock P.~204--206.

\bibitem{Malafeyev_rinz_corruption_agents_network_optimal_position}
\selectlanguageifdefined{english}
\BibEmph{Malafeev~O.~A., Redinskix~N.~D., Gerchiu~A.~L.} Optimizaczionnaya
  model' razmeshheniya korrupczionerov v seti~// Stroitel'stvo i
  e'kspluatacziya e'nergoe'ffektivny'x zdanij (teoriya i praktika s uchetom
  korrupczionnogo faktora)~/ L.~M.~Kolchedanczev, I.~N.~Legalov, G.~M.~Bad'in
  et~al. \BibDash
\newblock Borovichi~: NP ``NTO strojindustrii Sankt-Peterburga'', 2015.
  \BibDash
\newblock P.~128--140.

\bibitem{Malafeyev_rinz_compromise_auction_first_price_corruption}
\selectlanguageifdefined{english}
\BibEmph{Malafeev~O.~A., Koroleva~O.~A., Vasil'ev~YU.~G.} Kompromissnoe
  reshenie v aukczione pervoj czeny' s korrumpirovanny'm aukczionistom~//
  Stroitel'stvo i e'kspluatacziya e'nergoe'ffektivny'x zdanij (teoriya i
  praktika s uchetom korrupczionnogo faktora)~/ L.~M.~Kolchedanczev,
  I.~N.~Legalov, G.~M.~Bad'in et~al. \BibDash
\newblock Borovichi~: NP ``NTO strojindustrii Sankt-Peterburga'', 2015.
  \BibDash
\newblock P.~119--127.

\bibitem{Malafeyev_rinz_investment_with_possible_corruption}
\selectlanguageifdefined{english}
\BibEmph{Malafeev~O.~A., Redinskix~N.~D., Smirnova~T.~E.} Model'
  investirovaniya proekta s vozmozhnoj korrupcziej~// Stroitel'stvo i
  e'kspluatacziya e'nergoe'ffektivny'x zdanij (teoriya i praktika s uchetom
  korrupczionnogo faktora)~/ L.~M.~Kolchedanczev, I.~N.~Legalov, G.~M.~Bad'in
  et~al. \BibDash
\newblock Borovichi~: NP ``NTO strojindustrii Sankt-Peterburga'', 2015.
  \BibDash
\newblock P.~140--146.

\bibitem{Malafeyev_rinz_tender_real_estate_corruption}
\selectlanguageifdefined{english}
\BibEmph{Malafeev~O.~A., Axmady'shina~A.~R., Demidova~D.~A.} Model' tendera na
  ry'nke rie'lterskix uslug s uchetom korrupczii~// Stroitel'stvo i
  e'kspluatacziya e'nergoe'ffektivny'x zdanij (teoriya i praktika s uchetom
  korrupczionnogo faktora)~/ L.~M.~Kolchedanczev, I.~N.~Legalov, G.~M.~Bad'in
  et~al. \BibDash
\newblock Borovichi~: NP ``NTO strojindustrii Sankt-Peterburga'', 2015.
  \BibDash
\newblock P.~161--168.

\bibitem{Malafeyev_rinz_beginnings_geopolitics}
\selectlanguageifdefined{english}
\BibEmph{Malafeev~O.~A., Kefeli~I.~F.} Matematicheskie nachala global'noj
  geopolitiki. \BibDash
\newblock Sankt-Peterburg~: Sankt-Peterburgskij gosudarstvenny'j
  politexnicheskij universitet, 2013.

\bibitem{Malafeyev_rinz_evolution_mechanism_finnancial_economical_component}
\selectlanguageifdefined{english}
\BibEmph{Malafeev~O.~A., Maraxov~V.~G.} E'volyuczionny'j mexanizm dejstviya
  istochnikov i dvizhushhixsya sil grazhdanskogo obshhestva v sfere finansovoj
  i e'konomicheskoj komponenty' XXI veka~// K. Marks i budushhee filosofii
  Rossii~/ S.~V.~Busov, S.~I.~Dudnik, K.~YU.~ZHirkov et~al. \BibDash
\newblock Sankt-Peterburg~: OOO ``Izdatel'stvo VVM'', 2016. \BibDash
\newblock P.~112--135.

\bibitem{Malafeyev_rinz_stohastic_socio_economic_dynamics}
\selectlanguageifdefined{english}
\BibEmph{Malafeev~O.~A., Nemnyugin~S.~A.} Stoxasticheskaya model'
  soczial'no-e'konomicheskoj dinamiki~// Ustojchivost' i proczessy' upravleniya
  Materialy' III mezhdunarodnoj konferenczii. \BibDash
\newblock Sankt-Peterburg~: Izdatel'skij dom Fedorovoj G. V., 2015. \BibDash
\newblock P.~433--434.

\bibitem{Malafeyev_rinz_stohastic_estimation_firm_development_corruption}
\selectlanguageifdefined{english}
\BibEmph{Malafeev~O.~A., Redinskix~N.~D.} Stoxasticheskoe oczenivanie i prognoz
  e'ffektivnosti strategii razvitiya firmy' v usloviyax korrupczionnogo
  vozdejstviya~// Ustojchivost' i proczessy' upravleniya Materialy' III
  mezhdunarodnoj konferenczii. \BibDash
\newblock Sankt-Peterburg~: Izdatel'skij dom Fedorovoj G. V., 2015. \BibDash
\newblock P.~437--438.

\bibitem{Malafeyev_rinz_philosophical_beginnings_social_transformation}
\selectlanguageifdefined{english}
Filosofskie strategii soczial'ny'x preobrazovanij XXI veka~/ O.~A.~Malafeev,
  V.~N.~Volovich, T.~A.~Delieva et~al. \BibDash
\newblock Sankt-Peterburg~: Sankt-Peterburgskij gosudarstvenny'j universitet,
  2014.

\bibitem{Malafeyev_rinz_corrupt_company_federal_agency}
\selectlanguageifdefined{english}
\BibEmph{Malafeev~O.~A., Dejnega~L.~A., Andreeva~M.~A.} Model' vzaimodejstviya
  korrumpirovannogo predpriyatiya i federal'nogo otdela po bor'be s
  korrupcziej~// \BibEmph{Molodoj ucheny'j}. \BibDash
\newblock 2015. \BibDash
\newblock no. 12 (92). \BibDash
\newblock P.~15--20.

\bibitem{Malafeyev_rinz_dialog_mathematician_philosopher}
\selectlanguageifdefined{english}
\BibEmph{Malafeev~O.~A., Maraxov~V.~G.} Dialog filosofa i matematika: ``O
  filosofskix aspektax matematicheskogo modelirovaniya soczial'ny'x
  preobrazovanij XXI veka''~// Filosofiya poznaniya i tvorchestvo zhizni.
  \BibDash
\newblock Sankt-Peterburg~: Vladimir Dal', 2014. \BibDash
\newblock P.~279--292.

\bibitem{Malafeyev_rinz_modelling_customs}
\selectlanguageifdefined{english}
\BibEmph{Malafeev~O.~A., Drozdov~G.~D.} Modelirovanie tamozhennogo dela.
  \BibDash
\newblock Sankt-Peterburg~: Sankt-Peterburgskij gosudarstvenny'j universitet
  servisa i e'konomiki, 2013.

\bibitem{Malafeyev_rinz_global_arctic_game}
\selectlanguageifdefined{english}
\BibEmph{Malafeev~O.~A., Ivashov~L.~G., Kefeli~I.~F.} Global'naya arkticheskaya
  igra i ee uchastniki~// \BibEmph{Geopolitika i bezopasnost'}. \BibDash
\newblock 2014. \BibDash
\newblock no. 1 (25). \BibDash
\newblock P.~34--49.

\bibitem{Malafeyev_rinz_eurasian_arc_security_problem}
\selectlanguageifdefined{english}
Evrazijskaya duga nestabil'nosti i problemy' regional'noj bezopasnosti ot
  Vostochnoj Azii do Severnoj Afriki~/ O.~A.~Malafeev, M.~L.~Titarenko,
  L.~G.~Ivashov et~al. \BibDash
\newblock Sankt-Peterburg~: Studiya NP-Print, 2013.

\bibitem{Malafeyev_rinz_models_many_agents_geopolicy}
\selectlanguageifdefined{english}
\BibEmph{Malafeev~O.~A., Kefeli~I.~F.} O matematicheskix modelyax global'ny'x
  geopoliticheskix proczessov mnogoagentnogo vzaimodejstviya~//
  \BibEmph{Geopolitika i bezopasnost'}. \BibDash
\newblock 2013. \BibDash
\newblock no.~2. \BibDash
\newblock P.~44--57.

\bibitem{Malafeyev_rinz_defence_security}
\selectlanguageifdefined{english}
\BibEmph{Malafeev~O.~A., Kefeli~I.~F.} Nekotory'e zadachi obespecheniya
  oboronnoj bezopasnosti~// \BibEmph{Geopolitika i bezopasnost'}. \BibDash
\newblock 2013. \BibDash
\newblock no.~3. \BibDash
\newblock P.~84--92.

\bibitem{Malafeyev_rinz_linear_algebra_corruption}
\selectlanguageifdefined{english}
Linejnaya algebra s prilozheniyami k modelirovaniyu korrupczionny'x sistem i
  proczessov~/ O.~A.~Malafeev, N.N.~Sotnikova, I.V.~Zajczeva et~al. \BibDash
\newblock Stavropol'~: OOO ``Izdatel'skij dom ``TE'SE'RA'', 2016.

\bibitem{Malafeyev_rinz_control_dynamic_system}
\selectlanguageifdefined{english}
\BibEmph{Malafeev~O.~A.} Upravlenie v konfliktny'x dinamicheskix sistemax.
  \BibDash
\newblock Sankt-Peterburg~: SPbGU, 1993.

\bibitem{Malafeyev_scopus_2014_laser_radiation_control}
\selectlanguageifdefined{english}
Multi-criteria model of laser radiation control~/ O.A.~Malafeyev,
  E.G.~Neverova, S.A.~Nemnyugin, G.V.~Alferov~// 2014 2nd International
  Conference on Emission Electronics, ICEE 2014. \BibDash
\newblock New Jersey~: Institute of Electrical and Electronics Engineers Inc.,
  2014. \BibDash
\newblock P.~33--37.

\bibitem{Malafeyev_scopus_2014_postman_problem}
\selectlanguageifdefined{english}
\BibEmph{Malafeev~O., Grigorieva~X.} A competitive many-period postman problem
  with varying parameters~// \BibEmph{Applied Mathematical Sciences}. \BibDash
\newblock 2014. \BibDash
\newblock Vol.~8, no. 146. \BibDash
\newblock P.~7249--7258.

\bibitem{Malafeyev_scopus_2015_group_strategy_robots_interception}
\selectlanguageifdefined{english}
\BibEmph{Malafeyev~O., Alferov~G., Andreyeva~M.} Group strategy of robots in
  game-theoretic model of interception with incomplete information~// 2015
  International Conference on Mechanics - Seventh Polyakhov's Reading. \BibDash
\newblock New Jersey~: Institute of Electrical and Electronics Engineers Inc.,
  2015. \BibDash
\newblock P.~1--3.

\bibitem{Malafeyev_scopus_2015_programming_robots_inspection_interception}
\selectlanguageifdefined{english}
\BibEmph{Malafeyev~O.A., Alferov~G.V., Maltseva~A.S.} Programming the robot in
  tasks of inspection and interception~// 2015 International Conference on
  Mechanics - Seventh Polyakhov's Reading. \BibDash
\newblock New Jersey~: Institute of Electrical and Electronics Engineers Inc.,
  2015. \BibDash
\newblock P.~1--3.

\bibitem{Malafeyev_scopus_2014_charged_particle_beam}
\selectlanguageifdefined{english}
\BibEmph{Malafeyev~O.A., Nemnyugin~S.A., Alferov~G.V.} Charged particles beam
  focusing with uncontrollable changing parameters~// 2014 2nd International
  Conference on Emission Electronics, ICEE 2014. \BibDash
\newblock New Jersey~: Institute of Electrical and Electronics Engineers Inc.,
  2014. \BibDash
\newblock P.~25--27.

\bibitem{Malafeyev_scopus_2015_stochastic_socio_economic_model}
\selectlanguageifdefined{english}
\BibEmph{Malafeyev~O.A., Nemnyugin~S.A., Ivaniukovich~G.A.} Stochastic models
  of social-economic dynamics~// 2015 International Conference on ``Stability
  and Control Processes'' in Memory of V.I. Zubov, SCP 2015 - Proceedings.
  \BibDash
\newblock New Jersey~: Institute of Electrical and Electronics Engineers Inc.,
  2015. \BibDash
\newblock P.~483--485.

\bibitem{Malafeyev_scopus_2015_multicomponent_dynamics_single_sector_economy}
\selectlanguageifdefined{english}
\BibEmph{Malafeyev~O.A., Drozdov~G.D., Nemnyugin~S.A.} Multicomponent dynamics
  of competitive single-sector economy development~// 2015 International
  Conference on ``Stability and Control Processes'' in Memory of V.I. Zubov,
  SCP 2015 - Proceedings. \BibDash
\newblock New Jersey~: Institute of Electrical and Electronics Engineers Inc.,
  2015. \BibDash
\newblock P.~457--459.

\bibitem{Malafeyev_scopus_1996_robot_control_strategy}
\selectlanguageifdefined{english}
\BibEmph{Malafeyev~O.A., Alferov~G.V.} The robot control strategy in a domain
  with dynamical obstacles~// \BibEmph{Lecture Notes in Computer Science
  (including subseries Lecture Notes in Artificial Intelligence and Lecture
  Notes in Bioinformatics)}. \BibDash
\newblock 1996. \BibDash
\newblock Vol. 1093. \BibDash
\newblock P.~211--217.

\bibitem{Malafeyev_scopus_1996_system_in_external_field_stochastic}
\selectlanguageifdefined{english}
\BibEmph{Malafeev~O.A., Nemnyugin~S.A.} Generalized dynamic model of a system
  moving in an external field with stochastic components~//
  \BibEmph{Theoretical and Mathematical Physics}. \BibDash
\newblock 1996. \BibDash
\newblock Vol. 107, no.~3. \BibDash
\newblock P.~770--774.

\bibitem{Malafeyev_scopus_1974_electron_beam}
\selectlanguageifdefined{english}
EXPERIMENTAL STUDY OF AN ELECTRON BEAM IN DRIFT SPACE~/ M.A.~Vlasov,
  V.V.~Glebov, O.A.~Malafeyev, D.N.~Novichkov~// \BibEmph{Soviet journal of
  communications technology \& electronics}. \BibDash
\newblock 1986. \BibDash
\newblock Vol.~31, no.~3. \BibDash
\newblock P.~145--149.

\bibitem{Malafeyev_scopus_1974_existence_eps_equilibrium_dynamic_games}
\selectlanguageifdefined{english}
\BibEmph{Malafeev~O.A.} The existence of situations of yo-equilibrium in
  dynamic games with dependent movements~// \BibEmph{USSR Computational
  Mathematics and Mathematical Physics}. \BibDash
\newblock 1974. \BibDash
\newblock Vol.~14, no.~1. \BibDash
\newblock P.~88--99.

\bibitem{Malafeyev_scopus_1977_stationary_strategy}
\selectlanguageifdefined{english}
\BibEmph{Malafeev~O.A.} Stationary strategies in differential games~//
  \BibEmph{USSR Computational Mathematics and Mathematical Physics}. \BibDash
\newblock 1977. \BibDash
\newblock Vol.~17, no.~1. \BibDash
\newblock P.~37--46.

\bibitem{Malafeyev_scopus_1974_equilibrium_situations_dynamic_games}
\selectlanguageifdefined{english}
\BibEmph{Malafeev~O.A.} Equilibrium situations in dynamic games~//
  \BibEmph{Cybernetics}. \BibDash
\newblock 1974. \BibDash
\newblock Vol.~10, no.~3. \BibDash
\newblock P.~504--513.

\bibitem{Malafeyev_wos_2000_equilibrium_multicriteria_optimization}
\selectlanguageifdefined{english}
\BibEmph{Malafeyev~O.A., Troeva~M.S.} A weak equilibrium solution for
  multicriteria optimization problem~// Control Applications of Optimization
  2000: Proceedings of the 11th Ifac Workshop. \BibDash
\newblock Amsterdam~: Elsevier science \& technology, 2000. \BibDash
\newblock P.~363--368.

\bibitem{tumor_simos_icnaam}
\selectlanguageifdefined{english}
Numerical integration of Chaplain and stuart model~/ L.~Petrakis,
  Z.~Kalogiratou, Th.~Monovasilis, T.~E.~Simos~// \BibEmph{AIP Conference
  Proceedings}. \BibDash
\newblock 2016. \BibDash
\newblock Vol. 1738, no. 480131.

\bibitem{Kolpak2}
\selectlanguageifdefined{english}
\BibEmph{Kolpak~E.P., Ivanov~S.E.} On the three-dimensional Klein-Gordon
  equation with a cubic nonlinearity~// \BibEmph{International Journal of
  Mathematical Analysis}. \BibDash
\newblock 2016. \BibDash
\newblock Vol.~10, no. 13-16. \BibDash
\newblock P.~611--622.

\bibitem{Yeung2}
\selectlanguageifdefined{english}
\BibEmph{Yeung~D.W.K., Petrosyan~L.A.} Subgame consistent cooperative solution
  for NTU dynamic games via variable weights~// \BibEmph{Automatica}. \BibDash
\newblock 2015. \BibDash
\newblock Vol.~59. \BibDash
\newblock P.~84--89.

\end{thebibliography}
\bibliographystyle{ugost2008}%

\end{document}